\documentclass[letterpaper, 10 pt, conference]{ieeeconf}  

\usepackage{pdfsync}
\IEEEoverridecommandlockouts                              
\usepackage{graphics} 
\usepackage{epsfig} 
\usepackage{times} 
\usepackage{amsmath} 
\usepackage{amssymb}  

\newtheorem{thm}{Theorem}

\newtheorem{lemma}[thm]{Lemma}

\newtheorem{assum}{Assumption}

\newcommand{\R} {\mathcal{R}}

\newcommand{\w} {\mathbf{\omega}}

\newcommand{\bi} {\mathbf{b}_{\scriptscriptstyle i}}
\newcommand{\bj} {\mathbf{b}_{\scriptscriptstyle j}}
\newcommand{\bn}[1] {\mathbf{b}_{\scriptscriptstyle #1}}

\newcommand{\xone} {\mathbf{x_{\scriptscriptstyle 1}}}
\newcommand{\xtwo} {\mathbf{x_{\scriptscriptstyle 2}}}
\newcommand{\xn}[1] {x_{\scriptscriptstyle #1}}
\newcommand{\dxn}[1] {\dot{x}_{\scriptscriptstyle #1}}
\newcommand{\dxone} {\mathbf{\dot{x}_{\scriptscriptstyle 1}}}
\newcommand{\dxtwo} {\mathbf{\dot{x}_{\scriptscriptstyle 2}}}
\newcommand{\s} {\mathbf{s}}

\newcommand{\p} {\mathbf{p}}

\newcommand{\rpb} {\mathbf{r}}

\newcommand{\vel} {\mathbf{v}}
\newcommand{\velr} {\mathbf{v_{\scriptscriptstyle r}}}
\newcommand{\velc} {\mathbf{v_{\scriptscriptstyle c}}}
\newcommand{\dvelc} {\mathbf{\dot{v}_{\scriptscriptstyle c}}}

\newcommand{\rhoi} {\rho_{\scriptscriptstyle i}}
\newcommand{\drhoi} {\dot{\rho}_{\scriptscriptstyle i}}
\newcommand{\rhoj} {\rho_{\scriptscriptstyle j}}

\newcommand{\rhon}[1] {\rho_{\scriptscriptstyle #1}}

\newcommand{\skews}[1]{\mathcal{S}\left(#1\right)}

\DeclareMathOperator{\diag}{diag}

\newcommand{\eye}{\mathbf{I}}
\newcommand{\zeros}{\mathbf{0}}
\newcommand{\eyen}[1]{\mathbf{I}_{\scriptscriptstyle #1}}
\newcommand{\zerosn}[2]{\mathbf{0}_{\scriptscriptstyle #1 \times #2}}

\newcommand{\intu}{\mathbf{u}^{\scriptscriptstyle[1]}}
\newcommand{\intuT}{\mathbf{u}^{\scriptscriptstyle[1]\:\scriptstyle T}}

\newcommand{\ICF}{\{I\}\:}
\newcommand{\BCF}{\{B\}\:}

\title{Position USBL/DVL Sensor-based Navigation Filter in the presence of Unknown Ocean Currents}
\author{M. Morgado, P. Batista, P. Oliveira, and C. Silvestre%
\thanks{This work was partially supported by Funda\c{c}\~{a}o para a 
Ci\^{e}ncia e a Tecnologia (ISR/IST plurianual funding) and by the EU Project TRIDENT (Contract No. 
248497). The work of Marco Morgado was supported by PhD Student Scholarship SFRH/BD/25368/2005 from the Portuguese FCT POCTI programme.}
\thanks{The authors are with the Institute for Systems and Robotics, Instituto Superior T\'ecnico, Av. Rovisco Pais, 1, 1049-001, Lisbon, Portugal
{\tt\footnotesize \{marcomorgado,pbatista,pjcro,cjs\}@isr.ist.utl.pt}}}
\begin{document}
\maketitle
\thispagestyle{empty}
\pagestyle{empty}
\begin{abstract}
This paper presents a novel approach to the design of globally asymptotically stable (GAS) position filters for Autonomous Underwater Vehicles (AUVs) based directly on the nonlinear sensor readings of an Ultra-short Baseline (USBL) and a Doppler Velocity Log (DVL). Central to the proposed solution is the derivation of a linear time-varying (LTV) system that fully captures the dynamics of the nonlinear system, allowing for the use of powerful linear system analysis and filtering design tools that yield GAS filter error dynamics. Simulation results reveal that the proposed filter is able to achieve the same level of performance of more traditional solutions, such as the Extended Kalman Filter (EKF), while providing, at the same time, GAS guarantees, which are absent for the EKF.
\end{abstract}

\section{Introduction}
\label{sec:intro}

The design and implementation of navigation systems stands out as one of the most critical steps towards the successful operation of autonomous vehicles. The quality of the overall estimates of the navigation system dramatically influences the capability of the vehicles to perform precision-demanding tasks, see \cite{WhitcombICRA2000} and \cite{Oceans00ASIMOV} for interesting and detailed surveys on underwater vehicle navigation and its relevance. This paper presents a novel approach to the design of globally asymptotically stable (GAS) position filters directly based on the nonlinear sensor readings. 

Consider an underwater vehicle equipped with an Ultra-Short Baseline (USBL) underwater positioning device, a triad of orthogonally mounted rate gyros, and a Doppler Velocity Log (DVL), that moves in the presence of unknown ocean currents in a scenario that has a fixed transponder, as depicted in Fig. \ref{fig:MissionScenario_1transp2}. The USBL is composed of a small calibrated array of acoustic receivers, and measures the distance between the transponder and the receivers installed on-board. Given the proximity of the sensors in the receiving array, hence the name Ultra-Short Baseline (USBL), the USBL is capable of measuring more accurately the Range-Difference-of-Arrival (RDOA) of the acoustic waves at the receivers compared to the actual distances between the transponder and all the receivers. The DVL measures the velocity of the vehicle with respect to the fluid, and the rate gyros measure the angular velocities of the vehicle. Due to noisy measurements, unknown ocean currents, and the nonlinear nature of the range measurements, a filtering solution is required in order to correctly estimate the position of the transponder in the vehicle coordinate frame and the inertial velocity of the vehicle. Recent advances in the area of underwater navigation based on merging the information from acoustic arrays and other inertial sensors like DVLs, can be found in \cite{Rigby06}, \cite{Willemenot09}, \cite{BatistaAutomatica} and references therein.
\begin{center}
\vspace{-0.5cm}
\begin{figure}[t]
	\centering
		\includegraphics[width=0.58\columnwidth]{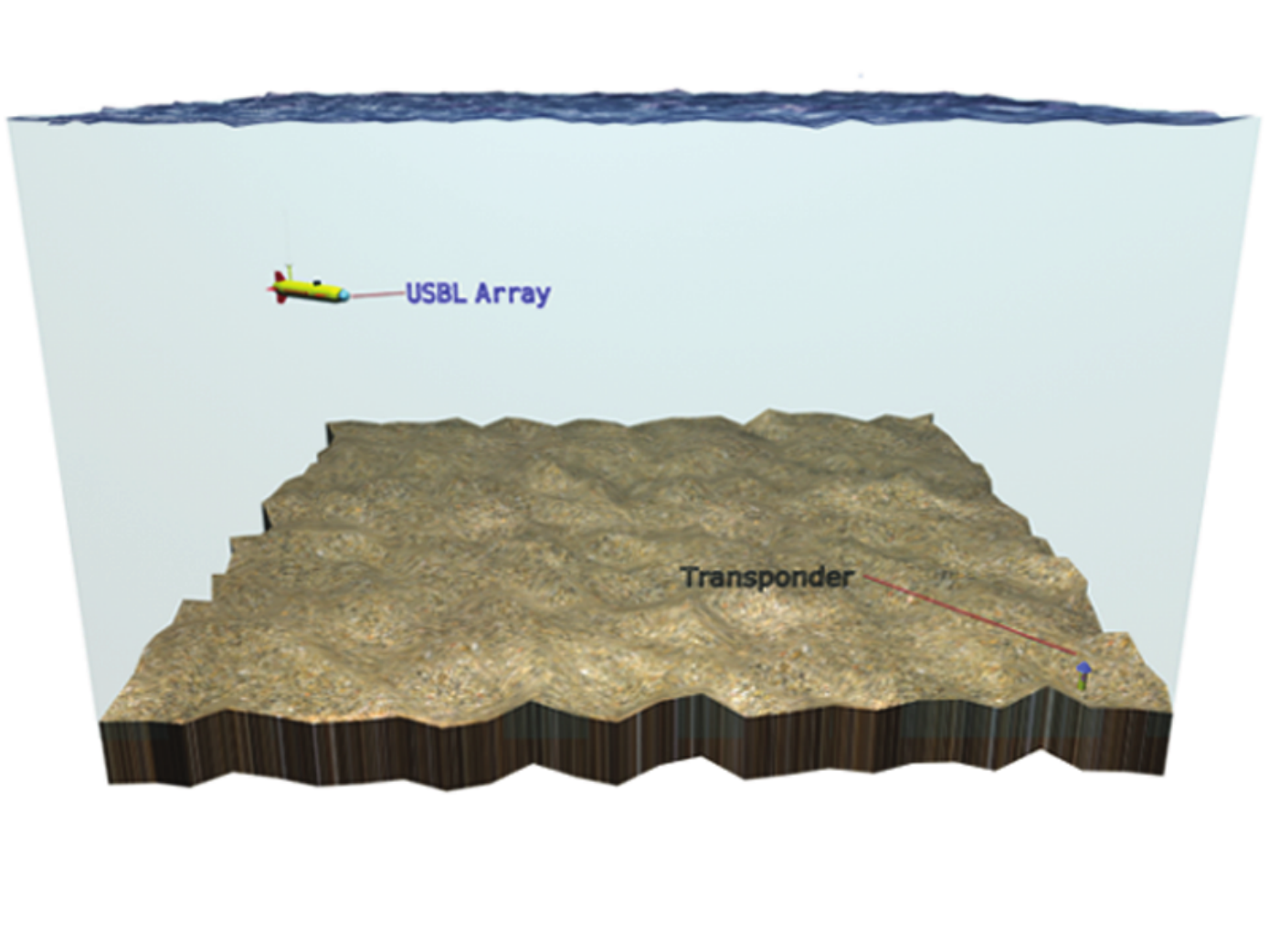}
		\vspace{-0.4cm}
	\caption{Mission scenario}
	\label{fig:MissionScenario_1transp2}
	\vspace{-1cm}
\end{figure}
\end{center}

Traditional solutions resort either to the well known Extended Kalman Filter (EKF)\cite{Fusion06}, Particle Filters (PF)\cite{Rigby06}, which lack global asymptotic stability properties, or to filtering solutions that use a precomputed position fix from the USBL device using the range and bearing and elevation angles of the transponder. The computation of this position fix often resorts to a Planar-Wave approximation of the acoustic wave arriving at the receiving array, previously used by the authors \cite{Fusion06}. In that case, the error cannot be guaranteed to converge to zero due to the planar-wave approximation. In fact, the error converges to a neighborhood of the origin, not arbitrarily small, that depends on the planar-wave approximation, and that only vanishes as the distance between the transponder and the vehicle approaches infinity. This behavior is obviously undesirable if the vehicle is to perform docking or other manoeuvres in the vicinity of a transponder, for instance.

The main contribution of this paper is the design of a globally asymptotically stable sensor-based filter to estimate the position of the transponder in the vehicle frame and the unknown current that biases the DVL readings. The solution presented in the paper departs from previous approaches as the range measurements are directly embedded in the filter structure, thus avoiding the planar-wave approximation, and follows related work found in \cite{batista2009single}, where single range measurements were considered and persistent excitation conditions were imposed on the vehicle motion to bear the system observable. In this paper the framework is extended to the case of having an array of receivers installed on-board the vehicle, which allows for the analysis of the overall system without any restriction on the vehicle motion. The nonlinear system dynamics are considered to their full extent and no linearizations are carried out whatsoever. At the core of the proposed filtering framework is the derivation of a linear time-varying (LTV) system that captures the dynamics of the nonlinear system. The LTV model is achieved through appropriate state augmentation, which is shown to mimic the nonlinear system, and ultimately allows for the use of powerful linear system analysis and filtering design tools that yield a novel estimation solution with GAS error dynamics. Although outside the scope of this work, the addition of an Attitude and Heading Reference System (AHRS) would allow for the final solution to be expressed directly in inertial coordinates.

The paper is organized as follows: Section \ref{sec:problemframework} sets the problem framework and definitions. The proposed filter design and main contributions of the paper are presented in Section \ref{sec:filter_design} where the filter structure is brought to full detail and an extensive observability analysis is carried out. Simulation results and performance comparison with traditional solutions are discussed in Section \ref{sec:simres}, and finally Section \ref{sec:conclusions} provides some concluding remarks.

\section{Problem framework}
\label{sec:problemframework}
In order to set the design framework, let $\ICF$ denote an inertial reference coordinate frame and $\BCF$ a coordinate frame attached to the vehicle, usually denominated as body-fixed coordinate frame. The position of the transponder $\rpb(t) \in \mathbb{R}^3$ in the vehicle coordinate frame $\BCF$ is given by
\begin{equation}
	\rpb(t) = \R^T(t)(\s - \p(t)),
\label{eq:bodyr1}
\end{equation} where $\s \in \mathbb{R}^3$ is the position of the transponder in inertial coordinates, $\p(t) \in \mathbb{R}^3$ is the position of the vehicle in inertial coordinates, and $\R(t) \in SO(3)$ is the rotation matrix from $\BCF$ to $\ICF$. The time derivative of $\R(t)$ verifies $\dot{\R}(t) = \R(t)\skews{\w(t)}$, where $\w(t) \in \mathbb{R}^3$ is the angular velocity of $\BCF$ with respect to $\ICF$ expressed in body-fixed coordinates, and $\skews{\w(t)}$ is the skew-symmetric matrix that represents the cross product such that $\skews{\w}a = \w \times a$.

Differentiating \eqref{eq:bodyr1} in time yields
\begin{align}
\dot{\rpb}(t)	&= -\skews{\w(t)}\rpb(t) - \vel(t), \label{eq:drdt}
\end{align} where $\vel(t) \in \mathbb{R}^3$ is the vehicle velocity expressed in body-fixed coordinates. The readings of the DVL are modeled by
\begin{equation}
\velr(t) = \vel(t)-\R^T(t){}^I\velc(t),
\label{eq:dvlreading}
\end{equation} where $\velr(t) \in \mathbb{R}^3$ is the velocity reading provided by the DVL, and ${}^I\velc(t) \in \mathbb{R}^3$ is the current velocity expressed in inertial coordinates and considered to be constant, that is, ${}^I\dvelc(t) = \zeros$. Using the current velocity expressed in body-fixed coordinates $\velc(t) = \R^T(t){}^I\velc(t)$ together with \eqref{eq:dvlreading} in \eqref{eq:drdt} yields
\begin{equation}
\dot{\rpb}(t)	= -\skews{\w(t)}\rpb(t) - \velc(t) - \velr(t).
\label{eq:drdt1}
\end{equation} The distances between the transponder and the receivers installed on-board the vehicle (as measured by the USBL) can be written as
\begin{align}
\rhoi(t) =& \|\bi - \rpb(t)\|,\:\:\:i=1,\ldots,n_r,\label{eq:rhoi1}
\end{align} where $\bi \in \mathbb{R}^3$ denotes the position of the receiver in $\BCF$, and $n_r$ is the number of receivers on the USBL. Combining the time-derivative of $\velc(t)$ with \eqref{eq:drdt1} and \eqref{eq:rhoi1} yields the nonlinear system
\begin{equation}\label{eq:nonlinearsys_body}
\left\{\begin{aligned}
\dot{\rpb}(t)	&= -\skews{\w(t)}\rpb(t) - \velc(t) - \velr(t), \\
\dvelc(t) &= -\skews{\w(t)}\velc(t), \\
\rhoi(t) &= \|\bi - \rpb(t)\|\mbox{, }i=1,\ldots,n_r.
\end{aligned}\right.
\end{equation}

The problem addressed in this paper is the design of a filter for \eqref{eq:nonlinearsys_body} considering noisy measurements.

\section{Filter design}
\label{sec:filter_design}

In this section the main results and contributions of the paper are presented. In order to reduce the complexity of the system dynamics a Lyapunov state transformation is firstly introduced in Section \ref{sec:state_transf}. The LTV system that will mimic the nonlinear behavior of the original system \eqref{eq:nonlinearsys_body} is proposed in Section \ref{sec:state_augmentation}, by means of an appropriate state augmentation. The observability analysis of the LTV system and its relation with the original nonlinear system is conducted in Section \ref{sec:obsanalysis}, and finally in Section \ref{sec:kalmanfilter}, the design of a Kalman filter is proposed in a stochastic setting for the resulting system.

\subsection{State transformation}
\label{sec:state_transf}

Consider the following state transformation
\begin{equation}
\begin{bmatrix} \xone(t) \\ \xtwo(t) \end{bmatrix}  := \mathbf{T}(t) \begin{bmatrix} \rpb(t) \\ \velc(t) \end{bmatrix},
\label{eq:statetransformation1}
\end{equation} where $\mathbf{T}(t) := \diag\left(\R(t),\R(t)\right) $ is a Lyapunov state transformation which preserves all observability properties of the original system \cite{brockett1970}.

The advantage of considering this state transformation is that the new unforced system dynamics becomes highly simplified as time-invariant, although the system output becomes time-varying and it is still nonlinear
\begin{equation}
\label{eq:nonlinearsys_inertial}
\left\{\begin{aligned}
\dxone(t)	&= - \xtwo(t) - \mathbf{u}(t), \\
\dxtwo(t) &= \zeros, \\
\rhoi(t) &= \|\bi - \R^T(t)\xone(t)\|,\:\:\: i=1,\ldots,n_r,
\end{aligned}\right.
\end{equation} where $\mathbf{u}(t) = \R(t)\velr(t)$.

\subsection{State augmentation}
\label{sec:state_augmentation}

In order to derive a linear system that mimics the dynamics of the original nonlinear system, a state augmentation procedure follows inherited directly from the kinematics of the nonlinear range outputs of \eqref{eq:nonlinearsys_inertial}. Thus, taking the time-derivative of $\rhoi(t)$ in \eqref{eq:nonlinearsys_inertial} yields
\begin{align}
&\drhoi(t) = \frac{1}{\rhoi(t)}\left[ (\bi^T\skews{\w(t)}\R^T(t) - \mathbf{u}^T(t))\xone(t) \right. \nonumber \\
&					 \left. + \bi^T\R^T(t)\xtwo(t) - \xone^T(t)\xtwo(t) +  \bi^T \R^T(t)\mathbf{u}(t)\right].
\label{eq:drhoi1}
\end{align}

Identifying the nonlinear part $\xone^T(t)\xtwo(t)$ in \eqref{eq:drhoi1} leads to the creation of the augmented states that will mimic this nonlinearity: $\xn{n_r + 3}(t) := \xone^T(t)\xtwo(t)$ and $\xn{n_r + 4}(t) := \|\xtwo(t)\|^2$, with the corresponding kinematics $\dxn{n_r + 3}(t) =  - \mathbf{u}^T(t)\xtwo(t) -\xn{n_r + 4}$ and $\dxn{n_r + 4}(t) = \zeros$.

Thus a new dynamic system is created by augmenting the original nonlinear system with the states
\begin{equation*}
\left\{
\begin{aligned}
\xn{3}(t) &:= \rhon{1}(t), \:\: \ldots \ldots ,&\xn{n_r + 2}(t) &:= \rhon{n_r}(t), \\
\xn{n_r + 3}(t) &:= \xone^T(t)\xtwo(t), &\xn{n_r + 4}(t) &:= \|\xtwo(t)\|^2,
\end{aligned}
\right.
\end{equation*} and denoting the new augmented state vector $\mathbf{x}(t) \in \mathbb{R}^{8+n_r}$ by
\begin{small}
\begin{align}
\mathbf{x}(t) &= \begin{bmatrix} \xone^T(t)\:\xtwo^T(t)\:\xn{3}(t)\:\ldots\:\xn{n_r + 2}(t)\:\xn{n_r + 3}(t)\:\xn{n_r + 4}(t) \end{bmatrix}^T. \nonumber
\end{align}
\end{small}

Combining the new augmented states dynamics with \eqref{eq:nonlinearsys_inertial} it is easy to verify that the augmented dynamics can be written as
\begin{equation*}
\mathbf{\dot{x}}(t) = \mathbf{A}(t)\mathbf{x}(t) + \mathbf{B}(t)\mathbf{u}(t),
\end{equation*} where
\begin{align}
\label{eq:sysdyn_matrixA}
&\mathbf{A}(t) = \nonumber \\
&\left[\begin{matrix}
\zeros & -\eye & \zeros & \zeros & \zeros \\
\zeros &  \zeros & \zeros & \zeros & \zeros \\
\frac{(\bn{1}^T\skews{\w(t)}\R^T(t) - \mathbf{u}^T(t))}{\rhon{1}(t)} &  \frac{\bn{1}^T\R^T(t)}{\rhon{1}(t)} & \zeros & -\frac{1}{\rhon{1}(t)} & 0 \\
\vdots & \vdots & \vdots & \vdots & \vdots \\
\frac{(\bn{n_r}^T\skews{\w(t)}\R^T(t) - \mathbf{u}^T(t))}{\rhon{n_r}(t)} &  \frac{\bn{n_r}^T\R^T(t)}{\rhon{n_r}(t)} & \zeros & -\frac{1}{\rhon{n_r}(t)} & 0 \\
\zeros & -\mathbf{u}^T(t) & \zeros & 0 & -1 \\
\zeros & \zeros & \zeros & 0 & 0
\end{matrix}
\right]
\end{align} and
\begin{equation}
\label{eq:sysdyn_matrixB}
\mathbf{B}(t) =
\left[\begin{matrix}
-\eyen{3} & \zeros & \frac{\R(t)\bn{1}}{\rhon{1}(t)} & \cdots & \frac{\R(t)\bn{n_r}}{\rhon{n_r}(t)} & \zeros & \zeros
\end{matrix}\right]^T.
\end{equation}

The following assumption is required so that \eqref{eq:sysdyn_matrixA} and \eqref{eq:sysdyn_matrixB} are well defined.

\begin{assum} \label{assum:vehicle_motion}
The motion of the vehicle is such that
\begin{displaymath}
\begin{array}{ccccc}
             \exists & \forall  & : & R_{min} \leq \rhoi(t) \leq R_{max},\\
             R_{min} > 0 & t \geq t_0 &  & \\
             R_{max} > 0 & i = 1,\ldots,n_r & & \\
\end{array}
\end{displaymath}
\end{assum}

From a practical point of view this is not restrictive since the vehicle and the coupled array will never be on top of a transponder, and neither will the ranges converge to infinity.

Note that the RDOA at the receivers are considered to be measured more accurately compared to the absolute distance between the transponder and any given reference receiver of the USBL. Selecting a reference sensor on the array, for instance receiver 1 for now, all the other ranges are easily reconstructed from the range measured at receiver 1 and the RDOA between receiver 1 and the other receivers, that is $\rhoj(t) = \rhon{1}(t) - \delta\rhon{1j}(t)$, where $\delta\rhon{1j}(t) = \rhon{1}(t) - \rhoj(t)$.

Taking into account that the augmented states $\xn{3}(t),\ldots,\xn{n_r + 2}(t)$ that correspond to the ranges are actually measured, it is straightforward to show from the outputs of \eqref{eq:nonlinearsys_body} that
\begin{align}
\rhoi^2(t) - \rhoj^2(t) = \|\bi\|^2 - \|\bj\|^2 - 2(\bi-\bj)^T\R^T(t)\xone(t), \nonumber
\end{align}
which leads to
\begin{align*}
\tfrac{2(\bi-\bj)^T\R^T(t)\xone(t)}{\rhoi(t) + \rhoj(t)}  + \rhoi(t) - \rhoj(t) = \tfrac{\|\bi\|^2 - \|\bj\|^2}{\rhoi(t) + \rhoj(t)} 
\end{align*} or, equivalently
\begin{equation}
\tfrac{2(\bi-\bj)^T\R^T(t)\xone(t)}{\rhoi(t) + \rhoj(t)}  + \xn{2+i}(t) - \xn{2+j}(t) = \tfrac{\|\bi\|^2 - \|\bj\|^2}{\rhoi(t) + \rhoj(t)},
\label{eq:squarediffequation1_equiv}
\end{equation} where the right hand-side of \eqref{eq:squarediffequation1_equiv} is measured and the left hand-side depends on the system state.

In order to complete the augmented system dynamics, discarding the original nonlinear outputs in \eqref{eq:nonlinearsys_inertial}, and considering \eqref{eq:squarediffequation1_equiv}, define the new augmented system outputs $\mathbf{y}(t) \in \mathbb{R}^{n_r+n_C}$ as
\begin{align*}
&\mathbf{y}(t) = \nonumber \\
&\begin{bmatrix}
\begin{bmatrix}
\xn{3}(t) &
\xn{3}(t) - \xn{4}(t) &
\cdots &
\xn{3}(t) - \xn{2+n_r}(t) \end{bmatrix}^T \\
\frac{2(\bn{1}-\bn{2})^T\R^T(t)\xone(t)}{\rhon{1}(t) + \rhon{2}(t)}  + \xn{2+1}(t) - \xn{2+2}(t) \\
\frac{2(\bn{1}-\bn{3})^T\R^T(t)\xone(t)}{\rhon{1}(t) + \rhon{3}(t)}  + \xn{2+1}(t) - \xn{2+3}(t) \\
\vdots\\
\frac{2(\bn{n_r-2}-\bn{n_r})^T\R^T(t)\xone(t)}{\rhon{n_r-2}(t) + \rhon{n_r}(t)} + \xn{2+n_r-2}(t) - \xn{2+n_r}(t) \\
\frac{2(\bn{n_r-1}-\bn{n_r})^T\R^T(t)\xone(t)}{\rhon{n_r-1}(t) + \rhon{n_r}(t)} + \xn{2+n_r-1}(t) - \xn{2+n_r}(t)
\end{bmatrix}, 
\end{align*} where $n_C = C_2^{n_r} = \frac{n_r!}{2(n_r - 2)!} =  \frac{n_r(n_r-1)}{2}$ is the number of all possible 2-combinations of $n_r$ elements.

In compact form, the augmented system dynamics can be written as
\begin{equation}
\left\{
\begin{aligned}
\mathbf{\dot{x}}(t) &= \mathbf{A}(t)\mathbf{x}(t) + \mathbf{B}(t)\mathbf{u}(t), \\
\mathbf{y}(t) &= \mathbf{C}(t)\mathbf{x}(t),
\end{aligned}
\right.
\label{eq:augmenteddynamics_fullmatrixform}
\end{equation} where
\begin{equation*}
\mathbf{C}(t) = 
\begin{bmatrix} 
\zerosn{n_r}{3} & \zerosn{n_r}{3} & \mathbf{C_0} & \zerosn{n_r}{2} \\ 
\mathbf{C_1}(t) & \zerosn{n_C}{3} & \mathbf{C_2} & \zerosn{n_C}{2}
\end{bmatrix}, 
\end{equation*}

\begin{equation*}
\mathbf{C_0} = \left[
\begin{smallmatrix}
1      & 0  & \cdots 	& 0 \\
1      & -1 &  				&  \\
\vdots &  \zeros  & \ddots  & \zeros  \\
1      &    &         & -1
\end{smallmatrix}\right],\:\:\:\:\:\:
\mathbf{C_2} = \left[
\begin{smallmatrix}
1      & -1     & 0   & 0      & \cdots 	 & 0  \\
1      & 0      & -1  & 0	     & \cdots	   & 0  \\
       &        &     & \vdots &           &    \\
0      & \cdots &  0  & 1      &  0        & -1 \\
0      & \cdots &  0  & 0      &  1        & -1
\end{smallmatrix}\right], 
\end{equation*} and
\begin{equation*}
\mathbf{C_1}(t) = 
\begin{bmatrix}
\frac{2(\bn{1}-\bn{2})^T\R^T(t)}{\rhon{1}(t) + \rhon{2}(t)} \\
\frac{2(\bn{1}-\bn{3})^T\R^T(t)}{\rhon{1}(t) + \rhon{3}(t)}\\
\vdots\\
\frac{2(\bn{n_r-2}-\bn{n_r})^T\R^T(t)}{\rhon{n_r-2}(t) + \rhon{n_r}(t)}\\
\frac{2(\bn{n_r-1}-\bn{n_r})^T\R^T(t)}{\rhon{n_r-1}(t) + \rhon{n_r}(t)}
\end{bmatrix}. 
\end{equation*}

\subsection{Observability analysis}
\label{sec:obsanalysis}

The Lyapunov state transformation and the state augmentation that were carried out allowed to derive the LTV system described in \eqref{eq:augmenteddynamics_fullmatrixform}, which ensembles the behavior of the original nonlinear system \eqref{eq:nonlinearsys_body}. The dynamic system \eqref{eq:augmenteddynamics_fullmatrixform} can be regarded as LTV, even though the system matrix $\mathbf{A}(t)$ depends explicitly on the system input and output, as evidenced by \eqref{eq:sysdyn_matrixA}. Nevertheless, this is not a problem from the theoretical point of view since both the input and output of the system are known continuous bounded signals. The idea is not new either, see, e.g., \cite{celikovsky2005secure}, and it just suggests, in this case, that the observability of \eqref{eq:augmenteddynamics_fullmatrixform} may be connected with the evolution of the system input or output (or both), which is not common and does not happen when this matrix does not depend on the system input or output. 

In order to fully understand and couple the behavior of both systems, the observability analysis of \eqref{eq:augmenteddynamics_fullmatrixform} is carried out in this section, using classical theory of linear systems. This analysis is conducted based on the observability Gramian associated with the pair $\left(\mathbf{A}(t),\mathbf{C}(t)\right)$, which is given by \cite{antsaklis2006linear}
\begin{equation*}
\mathcal{W}(t_0,t_f) = \int_{t_0}^{t_f}{\Phi^T(t,t_0)\mathbf{C}^T(t)\mathbf{C}(t)\Phi(t,t_0)}dt, 
\end{equation*} where $\Phi(t,t_0)$ is the state transition matrix of the LTV system \eqref{eq:augmenteddynamics_fullmatrixform}. Let
\begin{equation*}
\intu(t,t_0) = \int_{t_0}^{t}{\mathbf{u}(\sigma)d\sigma}. 
\end{equation*}

Tedious, lengthy, but straightforward, computations show that the transition matrix associated with $\mathbf{A}(t)$ is given by
\begin{align*}
\Phi(t,t_0) &= 
\begin{bmatrix}
\Phi_{AA}(t,t_0) & \zerosn{6}{n_r} 	&	\zerosn{6}{2} \\
\Phi_{BA}(t,t_0) & \eyen{n_r}			 	&	\Phi_{BC}(t,t_0) \\
\Phi_{CA}(t,t_0) & \zerosn{2}{n_r}	& \Phi_{CC}(t,t_0)
\end{bmatrix}, 
\end{align*} where
\begin{equation*}
\Phi_{AA}(t,t_0) = 
\begin{bmatrix}
\eye & -(t-t_0)\eye \\
\zeros & \eye
\end{bmatrix}, 
\end{equation*}
\begin{equation*}
\Phi_{BA}(t,t_0) = 
\begin{bmatrix}
\Phi_{BA1}(t,t_0) & \Phi_{BA2}(t,t_0)
\end{bmatrix},
\end{equation*}
\begin{equation*}
\Phi_{BA1}(t,t_0) = 
\begin{bmatrix}
\int_{t_0}^{t}{\frac{\bn{1}^T\skews{\w(\sigma)}\R^T(\sigma) - \mathbf{u}^T(\sigma)}{\rhon{1}(\sigma)}d\sigma} \\
\vdots \\
\int_{t_0}^{t}{\frac{\bn{n_r}^T\skews{\w(\sigma)}\R^T(\sigma) - \mathbf{u}^T(\sigma)}{\rhon{n_r}(\sigma)}d\sigma}
\end{bmatrix}, 
\end{equation*}
\begin{align*}
&\Phi_{BA2}(t,t_0) = \nonumber \\
&\begin{bmatrix}
\int_{t_0}^{t}{\frac{\bn{1}^T\R^T(\sigma)}{\rhon{1}(\sigma)}d\sigma} + \int_{t_0}^{t}{\frac{\intuT(\sigma,t_0)}{\rhon{1}(\sigma)}d\sigma} \\ + \int_{t_0}^{t}{\frac{-(\sigma-t_0)(\bn{1}^T\skews{\w(\sigma)}\R^T(\sigma) - \mathbf{u}^T(\sigma))}{\rhon{1}(\sigma)}d\sigma} \\
\vdots \\
\int_{t_0}^{t}{\frac{\bn{n_r}^T\R^T(\sigma)}{\rhon{n_r}(\sigma)}d\sigma} + \int_{t_0}^{t}{\frac{\intuT(\sigma,t_0)}{\rhon{n_r}(\sigma)}d\sigma} \\ + \int_{t_0}^{t}{\frac{-(\sigma-t_0)(\bn{n_r}^T\skews{\w(\sigma)}\R^T(\sigma) - \mathbf{u}^T(\sigma))}{\rhon{n_r}(\sigma)}d\sigma}
\end{bmatrix}, 
\end{align*}
\begin{align}
\Phi_{BC}(t,t_0) = 
\begin{bmatrix}
\Phi_{BC1}(t,t_0) & \Phi_{BC2}(t,t_0)
\end{bmatrix} = \nonumber \\
\begin{bmatrix}
-\int_{t_0}^{t}{\frac{1}{\rhon{1}(\sigma)}d\sigma} & \int_{t_0}^{t}{\frac{\sigma-t_0}{\rhon{1}(\sigma)}d\sigma} \\
\vdots & \vdots \\
-\int_{t_0}^{t}{\frac{1}{\rhon{n_r}(\sigma)}d\sigma} & \int_{t_0}^{t}{\frac{\sigma-t_0}{\rhon{n_r}(\sigma)}d\sigma}
\end{bmatrix}, \label{eq:transitionmatrix_1BCm}
\end{align} and $\Phi_{CA}(t,t_0)$ and $\Phi_{CC}(t,t_0)$ are omitted as they are not required in the sequel.

Before proceeding with the observability analysis, the following assumption is introduced which ultimately asserts the minimal number of receivers and configuration requirements of the USBL array in order to render the system observable regardless of the trajectory described of the vehicle.

\begin{assum}
\label{assum:receiver_configuration}
There are at least 4 non-coplanar receivers.
\end{assum}

The reasoning behind the need to have at least 4 non-coplanar receivers is that this is the minimal configuration that guarantees the uniqueness for the transponder position $\mathbf{r}(t)$.

The following theorem establishes the observability of the LTV system \eqref{eq:augmenteddynamics_fullmatrixform}.

\begin{thm}
\label{thm:observability_ltv}
The linear time-varying system \eqref{eq:augmenteddynamics_fullmatrixform} is observable on $[t_0,t_f]$, $t_0 < t_f$.

\begin{proof}
The observability proof of the LTV system \eqref{eq:augmenteddynamics_fullmatrixform} is accomplished by contradiction. Thus suppose that \eqref{eq:augmenteddynamics_fullmatrixform} is not observable on $\mathcal{I}:=[t_0,t_f]$. Then, there exists a non null vector $\mathbf{d} \in \mathbb{R}^{8+n_r}$
\begin{equation}
\mathbf{d} = \begin{bmatrix} \mathbf{d}_1^T & \mathbf{d}_2^T & \mathbf{d}_3^T & d_4 & d_5 \end{bmatrix},
\label{eq:definition_d}
\end{equation} with $\mathbf{d}_1 \in \mathbb{R}^{3}$, $\mathbf{d}_2 \in \mathbb{R}^{3}$, $\mathbf{d}_3 \in \mathbb{R}^{n_r}$, $d_4, d_5 \in \mathbb{R}$, such that $\mathbf{d}^T\mathcal{W}(t_0,t_f)\mathbf{d} = 0$ for all $t \in \mathcal{I}$, or equivalently,
\begin{eqnarray}
	\int_{t_0}^t{\| \mathbf{C}(\tau)\Phi(\tau,t_0) \mathbf{d}\|^2d\tau} = 0, & \forall_{t \in \mathcal{I}}.
	\label{eq:unobservability_condition1}
\end{eqnarray}
Taking the time derivative of \eqref{eq:unobservability_condition1} gives
\begin{eqnarray}
  \mathbf{C}(t)\Phi(t,t_0) \mathbf{d} = 0, & \forall_{t \in \mathcal{I}}.
	\label{eq:unobservability_condition2}
\end{eqnarray} From \eqref{eq:unobservability_condition2}, at $t=t_0$ comes
\begin{equation}
\begin{bmatrix} \mathbf{C}_0\mathbf{d}_3 \\ \mathbf{C}_1(t_0)\mathbf{d}_1 + \mathbf{C}_2\mathbf{d}_3 \end{bmatrix} = \zeros,
\label{eq:unobs_cond_part1}
\end{equation} which immediately implies that $\mathbf{C}_0\mathbf{d}_3 = \zeros$. As $\mathbf{C}_0$ is not singular allows to conclude that the only solution is the null vector
\begin{equation}
\mathbf{d}_3 = \zeros.
\label{eq:unobs_cond_sol_d3}
\end{equation} Replacing \eqref{eq:unobs_cond_sol_d3} in \eqref{eq:unobs_cond_part1} yields
\begin{equation}
\begin{bmatrix}
\frac{2}{\rhon{1}(t_0) + \rhon{2}(t_0)}(\bn{1}-\bn{2})^T \\
\frac{2}{\rhon{1}(t_0) + \rhon{3}(t_0)}(\bn{1}-\bn{3})^T \\
\vdots\\
\frac{2}{\rhon{n_r-2}(t_0) + \rhon{n_r}(t_0)}(\bn{n_r-2}-\bn{n_r})^T\\
\frac{2}{\rhon{n_r-1}(t_0) + \rhon{n_r}(t_0)}(\bn{n_r-1}-\bn{n_r})^T
\end{bmatrix}\R^T(t_0)\mathbf{d}_1 = \zeros.
\label{eq:unobs_cond_part2}
\end{equation} Under Assumption \ref{assum:receiver_configuration} the only solution for \eqref{eq:unobs_cond_part2} is
\begin{equation}
\mathbf{d}_1 = \zeros.
\label{eq:unobs_cond_sol_d1}
\end{equation} From \eqref{eq:unobservability_condition2} comes that
\begin{equation}
\mathbf{C_0}\Phi_{BA2}(t,t_0)\mathbf{d}_2 + \mathbf{C_0}\Phi_{BC}(t,t_0)\begin{bmatrix} d_4 \\ d_5 \end{bmatrix} = \zeros.
\label{eq:cphix_nR}
\end{equation} Taking the time derivative of \eqref{eq:cphix_nR} allows to write
\begin{align}
&\mathbf{C_0}\begin{bmatrix}
-\frac{(t-t_0)(\bn{1}^T\skews{\w(t)}\R^T(t) - \mathbf{u}^T(t))}{\rhon{1}(t)}\mathbf{d}_2 -\frac{d_4}{\rhon{1}(t)} + \\ \frac{\bn{1}^T\R^T(t)}{\rhon{1}(t)}\mathbf{d}_2 + \frac{\intuT(t,t_0)}{\rhon{1}(t)}\mathbf{d}_2  + \frac{(t-t_0)d_5}{\rhon{1}(t)} \\
\vdots\\
-\frac{(t-t_0)(\bn{n_r}^T\skews{\w(t)}\R^T(t) - \mathbf{u}^T(t))}{\rhon{n_r}(t)}\mathbf{d}_2 -\frac{d_4}{\rhon{n_r}(t)} + \\ \frac{\bn{n_r}^T\R^T(t)}{\rhon{n_r}(t)}\mathbf{d}_2 +  \frac{\intuT(t,t_0)}{\rhon{n_r}(t)}\mathbf{d}_2  + \frac{(t-t_0)d_5}{\rhon{n_r}(t)}
\end{bmatrix} = \zeros.
\label{eq:dotcphix_nR}
\end{align} Evaluating \eqref{eq:dotcphix_nR} at $t=t_0$ yields
\begin{align}
\mathbf{C_0}\begin{bmatrix}
\frac{\bn{1}^T}{\rhon{1}(t_0)} & \frac{-1}{\rhon{1}(t_0)} \\
\vdots & \vdots \\
\frac{\bn{n_r}^T}{\rhon{n_r}(t_0)} & \frac{-1}{\rhon{n_r}(t_0)}
\end{bmatrix}\begin{bmatrix} \R^T(t_0)\mathbf{d}_2 \\ d_4 \end{bmatrix} = \zeros.
\label{eq:dotcphix_nR2}
\end{align} Taking into account that $\mathbf{C}_0$ is not singular, it is easy to verify that under Assumption \ref{assum:receiver_configuration} the only solution for \eqref{eq:dotcphix_nR2} is
\begin{equation}
\begin{bmatrix}
\frac{\bn{1}^T}{\rhon{1}(t_0)} & \frac{-1}{\rhon{1}(t_0)} \\
\vdots & \vdots \\
\frac{\bn{n_r}^T}{\rhon{n_r}(t_0)} & \frac{-1}{\rhon{n_r}(t_0)}
\end{bmatrix}\begin{bmatrix} \R^T(t_0)\mathbf{d}_2 \\ d_4 \end{bmatrix} = \zeros \Rightarrow \begin{bmatrix} \mathbf{d}_2 \\ d_4 \end{bmatrix} = \zeros.
\label{eq:unobs_cond_sol_d2d4}
\end{equation} Now setting \eqref{eq:unobs_cond_sol_d3}, \eqref{eq:unobs_cond_sol_d1}, and \eqref{eq:unobs_cond_sol_d2d4} in \eqref{eq:dotcphix_nR} yields
\begin{align}
\mathbf{C_0}\begin{bmatrix}
\frac{t-t_0}{\rhon{1}(t_0)}d_5 \\
\vdots\\
\frac{t-t_0}{\rhon{n_r}(t_0)}d_5
\end{bmatrix} = \zeros.
\label{eq:dotcphix_nR3}
\end{align} Again the only possible solution for \eqref{eq:dotcphix_nR3} is
\begin{equation*}
d_5 = 0. 
\end{equation*} This concludes the proof since the only solution $\mathbf{d} = \zeros$ of \eqref{eq:unobservability_condition1} contradicts the hypothesis of the existence of a non null vector $\mathbf{d}$ such that \eqref{eq:unobservability_condition1} is true. Thus, by contradiction, the LTV system \eqref{eq:augmenteddynamics_fullmatrixform} is observable.
\end{proof}
\end{thm}

Although the observability of the LTV system \eqref{eq:augmenteddynamics_fullmatrixform} has been established, it does not mean that the original nonlinear system \eqref{eq:nonlinearsys_body} is also observable, and neither means that an observer for \eqref{eq:augmenteddynamics_fullmatrixform} is also an observer for \eqref{eq:nonlinearsys_body}. This however turns out to be true, as it is shown in the next theorem.

\begin{thm} \label{thm:nonlinear_obs}
The nonlinear system \eqref{eq:nonlinearsys_inertial} is observable in the sense that, given $\{\mathbf{y}(t),t\in[t_0,t_f]\}$ and $\{\mathbf{u}(t),t\in[t_0,t_f]\}$, the initial state $\mathbf{x}(t_0) = \begin{bmatrix} \mathbf{x}_1^T(t_0) & \mathbf{x}_2^T(t_0) \end{bmatrix}^T$ is uniquely defined. Moreover, a state observer for the LTV system \eqref{eq:augmenteddynamics_fullmatrixform} with globally asymptotically stable error dynamics is also a state observer for the nonlinear system \eqref{eq:nonlinearsys_inertial}, with globally asymptotically stable error dynamics.

\begin{proof}
The observability of the LTV system \eqref{eq:augmenteddynamics_fullmatrixform} has already been established in Theorem \ref{thm:observability_ltv}, thus given $\{\mathbf{y}(t),t\in[t_0,t_f]\}$ and $\{\mathbf{u}(t),t\in[t_0,t_f]\}$, the initial state of \eqref{eq:augmenteddynamics_fullmatrixform} is uniquely defined. Let
\begin{equation*}
\mathbf{z}(t_0) = \begin{bmatrix} \mathbf{z}_1^T(t_0) \: \mathbf{z}_2^T(t_0) \: \mathbf{z}_3^T(t_0) \: z_4(t_0) \: z_5(t_0) \end{bmatrix}^T, 
\end{equation*} with $\mathbf{z}_1(t_0),\mathbf{z}_2(t_0)\in \mathbb{R}^3$, $\mathbf{z}_3(t_0)\in \mathbb{R}^{n_r}$, and $z_4(t_0), z_5(t_0)\in \mathbb{R}^3$ be the initial state of the LTV system \eqref{eq:augmenteddynamics_fullmatrixform} and 
\begin{equation*}
\mathbf{x}(t_0) = \begin{bmatrix} \mathbf{x}_1^T(t_0) \: \mathbf{x}_2^T(t_0) \end{bmatrix}^T, 
\end{equation*} be the initial state of the nonlinear system \eqref{eq:nonlinearsys_inertial}. The evolution of $\mathbf{x}_1(t)$ for the nonlinear system \eqref{eq:nonlinearsys_inertial} can be easily shown to be given by
\begin{equation}
\mathbf{x}_1(t) = \mathbf{x}_1(t_0) - (t-t_0)\mathbf{x}_2(t_0) - \intu(t,t_0),
\label{eq:nlsys_evolx1}
\end{equation} which is similar to the evolution of $\mathbf{x}_1(t)$ for the LTV system differing only in the initial condition. Using \eqref{eq:nlsys_evolx1}, the output of the nonlinear system \eqref{eq:nonlinearsys_inertial} can be shown to satisfy
\begin{align}
\rhoi^2(t) =& \|\mathbf{x}_1(t_0) - \R(t_0)\bi\|^2 - 2\mathbf{x}_1^T(t_0)\R(t)\bi \nonumber \\ 
 &+ (t-t_0)^2\|\mathbf{x}_2(t_0)\|^2 + 2\mathbf{x}_1^T(t_0)\R(t_0))\bi \nonumber \\
 &-2(t-t_0)\mathbf{x}_2^T(t_0)\mathbf{x}_1(t_0) + 2(t-t_0)\mathbf{x}_2^T(t_0)\R(t)\bi \nonumber \\
 &+2(t-t_0)\mathbf{x}_2^T(t_0)\intu(t,t_0) - 2\intuT(t,t_0)\mathbf{x}_1(t_0) \nonumber \\
 &+2\intuT(t,t_0)\R(t)\bi + \|\intu(t,t_0)\|^2,
\label{eq:rhoi2nl_2}
\end{align} and the squared range difference between receiver $i$ and $j$
\begin{align}
&\rhoi^2(t) - \rhoj^2(t) = 2(\bi - \bj)^T(\R^T(t_0)\mathbf{x}_1(t_0) - \R^T(t)\mathbf{x}_1(t)) \nonumber \\
 &  + |\mathbf{x}_1(t_0) - \R(t_0)\bi\|^2 - \|\mathbf{x}_1(t_0) - \R(t_0)\bj\|^2.
\label{eq:rhoij2nl}
\end{align} Now, multiplying the set of augmented outputs of the LTV system \eqref{eq:augmenteddynamics_fullmatrixform} by the corresponding sum of pair of ranges, and taking into account that the states $\xn{3}(t),\ldots,\xn{2+n_r}(t)$ correspond to the actual ranges $\rhon{1}(t),\ldots,\rhon{n_r}(t)$, gives that
\begin{align}
&\left[\tfrac{2(\bi-\bj)^T\R^T(t)\xone(t)}{\rhoi(t) + \rhoj(t)}  + \xn{2+i}(t) - \xn{2+j}(t)\right]\left(\rhoi(t) + \rhoj(t)\right) \nonumber \\
& = 2(\bi-\bj)^T\R^T(t)\xone(t)  + \rhoi^2(t) - \rhoj^2(t).
\label{eq:augmented_setmult}
\end{align} The squared range of the LTV system \eqref{eq:augmenteddynamics_fullmatrixform} in \eqref{eq:augmented_setmult} can be shown to satisfy
\begin{align}
\rhoi^2(t) =& z_{2+i}^2(t_0) - 2\mathbf{z}_1^T(t_0)\R(t)\bi \nonumber \\
 &+(t-t_0)^2z_{n_r+4}(t_0)  + 2\mathbf{z}_1^T(t_0)\R(t_0)\bi \nonumber \\
 &-2(t-t_0)z_{n_r+3}(t_0) + 2(t-t_0)\mathbf{z}_2^T(t_0)\R(t)\bi \nonumber \\
 &+2(t-t_0)\mathbf{z}_2^T(t_0)\intu(t,t_0) - 2\intuT(t,t_0)\mathbf{z}_1(t_0) \nonumber \\
 &+2\intuT(t,t_0)\R(t)\bi + \|\intu(t,t_0)\|^2,
\label{eq:rhoi2ltv_1}
\end{align} and consequently it is true, for the LTV system \eqref{eq:augmenteddynamics_fullmatrixform}, that
\begin{align}
&\rhoi^2(t) - \rhoj^2(t) = z_{2+i}^2(t_0) - z_{2+j}^2(t_0)  \nonumber \\ 
&- 2(\bi-\bj)^T\R^T(t)\left(\mathbf{z}_1(t_0) - (t-t_0)\mathbf{z}_2^T(t_0) - \intu(t,t_0)\right) \nonumber \\
& + 2(\bi-\bj)^T\R^T(t_0)\mathbf{z}_1(t_0). \label{eq:rhoij2ltv}
\end{align} From the LTV system \eqref{eq:augmenteddynamics_fullmatrixform} transition matrix and forced response comes
\begin{equation}
\mathbf{x}_1(t) = \mathbf{z}_1(t_0) - (t-t_0)\mathbf{z}_2(t_0) - \intu(t,t_0).
\label{eq:evol_x1_ltv}
\end{equation} Replacing \eqref{eq:evol_x1_ltv} and \eqref{eq:rhoij2ltv} in \eqref{eq:augmented_setmult} yields
\begin{align}
2(\bi-\bj)^T\R^T(t)\xone(t)  + \rhoi^2(t) - \rhoj^2(t) = \nonumber \\
z_{2+i}^2(t_0) - z_{2+j}^2(t_0) + 2\mathbf{z}_1^T(t_0)\R(t_0)(\bi-\bj).
\label{eq:augmented_setmult_evol}
\end{align} The augmented states of the LTV system \eqref{eq:augmenteddynamics_fullmatrixform} $x_{2+i}(t)$, $i = 1,\ldots,n_r$, are actually considered to be measured and correspond to the range measurements. Therefore, it must be for its initial condition
\begin{equation}
z_{2+i}(t_0) = \|\mathbf{x}_1(t_0) - \R(t_0)\bi\|.
\label{eq:init_cond_ltvaugmented}
\end{equation} Now taking \eqref{eq:init_cond_ltvaugmented} into account, and comparing the evolution of the augmented outputs \eqref{eq:augmented_setmult_evol} for the LTV system with \eqref{eq:rhoij2nl} for the nonlinear system, it follows that
\begin{equation}
\mathbf{B_c}\R^T(t_0)\left[ \mathbf{x}_1(t_0) - \mathbf{z}_1(t_0) \right] = \zeros,
\label{eq:obs_ltv2nl_cond1}
\end{equation} where
\begin{equation*}
\mathbf{B_c} = \begin{bmatrix}
(\bn{1}-\bn{2})^T \\
(\bn{1}-\bn{3})^T \\
\vdots\\
(\bn{n_r-2}-\bn{n_r})^T \\
(\bn{n_r-1}-\bn{n_r})^T \\
\end{bmatrix}. 
\end{equation*} Under Assumption \ref{assum:receiver_configuration} the only solution of \eqref{eq:obs_ltv2nl_cond1} is
\begin{equation}
\mathbf{x}_1(t_0) = \mathbf{z}_1(t_0).
\label{eq:obs_ltv2nl_sol1}
\end{equation} Setting \eqref{eq:init_cond_ltvaugmented} and \eqref{eq:obs_ltv2nl_sol1} and comparing the difference between square of ranges for both systems yields
\begin{equation}
2(t-t_0)(\bi-\bj)^T\R^T(t)\left[ \mathbf{x}_2(t_0) - \mathbf{z}_2(t_0) \right] = \zeros.
\label{eq:obs_ltv2nl_cond2}
\end{equation} Taking the time derivative of \eqref{eq:obs_ltv2nl_cond2} comes
\begin{align}
	&\left[- 2(t-t_0)(\bi-\bj)^T\skews{\w(t)}\R^T(t) \right. \nonumber \\ &\left.+ 2(\bi-\bj)^T\R^T(t)\right]\left( \mathbf{x}_2(t_0) - \mathbf{z}_2(t_0) \right) = \zeros. \label{eq:obs_ltv2nl_cond3}
\end{align} At $t=t_0$, \eqref{eq:obs_ltv2nl_cond3} comes as
\begin{equation}
\mathbf{B_c}\R^T(t_0)\left[ \mathbf{x}_2(t_0) - \mathbf{z}_2(t_0) \right] = \zeros.
\label{eq:obs_ltv2nl_cond4}
\end{equation} Again under Assumption \ref{assum:receiver_configuration} the only solution of \eqref{eq:obs_ltv2nl_cond4} is
\begin{equation}
\mathbf{x}_2(t_0) = \mathbf{z}_2(t_0).
\label{eq:obs_ltv2nl_sol2}
\end{equation} Finally, setting \eqref{eq:init_cond_ltvaugmented}, \eqref{eq:obs_ltv2nl_sol1}, and \eqref{eq:obs_ltv2nl_sol2} in \eqref{eq:rhoi2ltv_1} and comparing to \eqref{eq:rhoi2nl_2} yields
\begin{align}
&-2(t-t_0)\left(\mathbf{x}_2^T(t_0)\mathbf{x}_1(t_0) - z_{n_r+3}(t_0)\right) \nonumber \\
&+(t-t_0)^2\left(\|\mathbf{x}_2(t_0)\|^2 - z_{n_r+4}(t_0)\right) = \zeros.
\label{eq:obs_ltv2nl_cond5}
\end{align} As $(t-t_0)$ and $(t-t_0)^2$ are linearly independent functions, the only solution for \eqref{eq:obs_ltv2nl_cond5} is $z_{n_r+3}(t_0) = \mathbf{x}_2^T(t_0)\mathbf{x}_1(t_0)$ and $z_{n_r+4}(t_0) = \|\mathbf{x}_2(t_0)\|^2$. Thus, the initial state of the nonlinear system \eqref{eq:nonlinearsys_inertial} matches the initial state of the LTV system \eqref{eq:augmenteddynamics_fullmatrixform} which is uniquely defined. Therefore the nonlinear system \eqref{eq:nonlinearsys_inertial} is also observable. 
\end{proof}
\end{thm}

Note that the usual concept of observability for nonlinear systems is not as strong as that presented in the statement of Theorem \ref{thm:nonlinear_obs}, see \cite{Hermann1977}. Also, in the observer structure that was derived, there is nothing imposing the nonlinear algebraic relations between the original and the additional states, allowing for the dynamics of the new system to be considered linear time-varying. It is however trivial to show that those relations are indeed preserved asymptotically.

Although the observability results were derived with respect to the nonlinear system \eqref{eq:nonlinearsys_inertial}, they also apply to the original nonlinear system \eqref{eq:nonlinearsys_body} as they are related through a Lyapunov transformation. Thus, the design of an observer for the original nonlinear system follows simply by reversing the state transformation \eqref{eq:statetransformation1}, as it will be detailed in the following section.

\subsection{Kalman filter}
\label{sec:kalmanfilter}

The observer structure devised so far was based on a deterministic setting providing strong constructive results, in the sense that it was shown, in Theorem \ref{thm:nonlinear_obs}, that an observer with globally asymptotically stable error dynamics for the LTV system \eqref{eq:augmenteddynamics_fullmatrixform} provides globally asymptotically stable error dynamics for the estimation of the state of the original nonlinear system. However, in practice there exists measurement noise and system disturbances, motivating the derivation of a filtering solution within a stochastic setting. Therefore, the design of a LTV Kalman Filter (even tough other filtering solutions could be used, e.g. a $\mathcal{H}_\infty$ filter) is presented next. Before proceeding with the derivation of the proposed filter, it is important to stress, however, that this filter is not optimal, as the existence of multiplicative noise is evident by looking into the LTV system matrices.

Nevertheless, the errors associated with the Kalman filter estimates are GAS, as it can be shown that the system is not only observable but also uniformly completely observable, a sufficient condition for the stability of the LTV Kalman filter \cite{Anderson1971}. The following assumption and lemma are introduced to guarantee the uniform complete observability of the system.

\begin{assum}\label{assum:bounded_signals}
The position of the transponder in the vehicle coordinate frame $\rpb(t)$, and the angular and linear velocities, $\w(t)$ and $\vel(t)$ respectively, are bounded signals. Moreover, the time derivatives of these signals ($\dot{\rpb}(t)$, $\dot{\w}(t)$, and $\dot{\vel}(t)$ respectively), are also bounded, as well as the derivatives of the ranges $\drhoi(t)$, with $i=1,\ldots,n_r$.
\end{assum}

\begin{lemma}[{\cite[Lemma 1]{paper:Batista:CDC:2009}}] \label{lemma:batista_intfunc}
     Let
     $\mathbf{f}(t) : \left[t_0, t_f\right] \subset \mathbb{R} \to \mathbb{R}^n$
     be a continuous and two times continuously differentiable function
     on     ${\cal I} := \left[t_0, t_f\right]$, $T := t_f - t_0 > 0$, and
     such that $\mathbf{f}\left(t_0\right) =  \mathbf{0}$. Further assume that
     $\max_{t \in {\cal I}} \| \ddot{\mathbf{f}}(t)\| \leq C$.
     If there exists a constant $\alpha^* > 0$ and a time $t^* \in {\cal I}$
     such that $\| \dot{\mathbf{f}}(t^*) \| \geq \alpha^*$, then there exists 
     constants $\beta^* > 0$ and $0 < \delta^* \leq T$ such that 
     $\| \mathbf{f}\left(t_0 +\delta^* \right) \| \geq \beta^*$.
     \label{thm:ODSA:Lemma}
\end{lemma}

From a practical point of view, Assumption \ref{assum:bounded_signals} is not restrictive since the systems presented herein are in fact finite energy systems that ensemble realizable physical vehicles and sensors. The LTV system \eqref{eq:augmenteddynamics_fullmatrixform} is finally shown to be uniformly completely observable in the following theorem.

\begin{thm}
The linear time-varying system \eqref{eq:augmenteddynamics_fullmatrixform} is uniformly completely observable, that is, there exists positive constants $\alpha_1$, $\alpha_2$, $\delta$, such that $\alpha_1\eye \preceq \mathcal{W}(t,t+\delta) \preceq \alpha_2\eye$ for all $t \geq t_0$.
\begin{proof} The bounds on the observability Gramian $\mathcal{W}(t,t+\delta)$ can be written as
\begin{equation}
\alpha_1 \leq \mathbf{d}^{\scriptscriptstyle T}\mathcal{W}(t,t+\delta)\mathbf{d} \leq \alpha_2,
\label{eq:bound_gramian_uniform}
\end{equation} for all $t \geq t_0$, and for all $\mathbf{d} \in \mathbb{R}^{8+n_r}$ such that $\|\mathbf{d}\|=1$. The proof follows by noticing that \eqref{eq:bound_gramian_uniform} can be written as
$\alpha_1 \leq \int_{t}^{t+\delta}{\|\mathbf{f}(\tau)\|^2d\tau} \leq \alpha_2$, where
\begin{equation}
\mathbf{f}(\tau) := \mathbf{C}(\tau)\mathbf{\Phi}(\tau,t)\mathbf{d}.
\label{eq:uniform_f1}
\end{equation} The existence of the upper bound $\alpha_2$ is trivially checked, as under Assumption \ref{assum:bounded_signals} the matrices $\mathbf{A}(t)$ and $\mathbf{C}(t)$ are norm-bounded and $\mathbf{f}(\tau)$ is integrated over limited intervals. Let $\mathbf{d} = \begin{bmatrix} \mathbf{d}_1^{\scriptscriptstyle T} & \mathbf{d}_2^{\scriptscriptstyle T} & \mathbf{d}_3^{\scriptscriptstyle T} & d_4 & d_5 \end{bmatrix}$, with $\mathbf{d}_1 \in \mathbb{R}^{3}$, $\mathbf{d}_2 \in \mathbb{R}^{3}$, $\mathbf{d}_3 \in \mathbb{R}^{n_r}$, $d_4, d_5 \in \mathbb{R}$. Evaluating \eqref{eq:uniform_f1} at $\tau=t$, it is straightforward to verify that if $\mathbf{d}_3\neq\zeros$, then $\|\mathbf{f}(t)\|$ is immediately bounded. Indeed, as $\|\mathbf{f}(t)\| \geq \|\mathbf{C_0}\mathbf{d}_3\| = \alpha_1^* > 0$  for all $t \geq t_0$, since $\mathbf{C_0}$ has full column rank by construction. Suppose now that $\mathbf{d}_3 = \zeros$. Then, it can also be seen that if $\mathbf{d}_1\neq\zeros$, it is true that $\|\mathbf{f}(t)\| = \|\mathbf{C_1}(t)\mathbf{d}_1\|$, which is clearly bounded by
\begin{equation}
\|\mathbf{f}(t)\| \geq \sigma_{min}(\mathbf{D}_{\rho+}^{-1}(t)) \|2 \mathbf{C_2} \mathbf{U_r} \R^{\scriptscriptstyle T}(t)\mathbf{d}_1\|,
\end{equation}  for all $t \geq t_0$, where the operator $\sigma_{min}(\mathbf{A})$ represents the smallest singular value of $\mathbf{A}$ and
\begin{align*}
\mathbf{U_r} &:= \begin{bmatrix} \bn{1} & \cdots & \bn{n_r}\end{bmatrix}^{\scriptscriptstyle T}  & \in \mathbb{R}^{n_r \times 3}, \\
\mathbf{D}_{\rho+}(t) &:=  \diag\left(\left[
\begin{smallmatrix}
\rhon{1}(t) + \rhon{2}(t) \\
\rhon{1}(t) + \rhon{3}(t)\\
\vdots\\
\rhon{n_r-2}(t) + \rhon{n_r}(t)\\
\rhon{n_r-1}(t) + \rhon{n_r}(t)
\end{smallmatrix}\right]\right) & \in \mathbb{R}^{n_c \times n_c}.
\end{align*}
Under Assumption \ref{assum:vehicle_motion} it is clear that $\sigma_{min}(\mathbf{D}_{\rho+}^{-1}(t)) \geq \frac{1}{R_{max}}$ for all $t \geq t_0$, which implies
\begin{align*}
\|\mathbf{f}(t)\| \geq \frac{1}{R_{max}}\sigma_{min}(\mathbf{C_2} \mathbf{U_r})\left\|\R^T(t)\mathbf{d}_1\right\|
\end{align*}  for all $t \geq t_0$. Now under Assumption \ref{assum:receiver_configuration} it follows that $\mathbf{C_2} \mathbf{U_r}$ has full column rank and therefore there exists a positive constant $\beta_1^*$ such that $\sigma_{min}(\mathbf{C_2} \mathbf{U_r}) = \beta_1^* > 0$. Thus, taking into account that $\left\|\R^T(t)\mathbf{d}_1\right\| = \left\|\mathbf{d}_1\right\|$,
\begin{align}
\|\mathbf{f}(t)\| &\geq \frac{\beta_1^*}{R_{max}}\left\|\mathbf{d}_1\right\| = \alpha_2^* > 0 
\end{align} for all $t \geq t_0$. Using the same set of assumptions and procedures it is possible to show that the derivative of $\mathbf{f}(\tau)$ evaluated at $\tau = t$ is also uniformly bounded
\begin{align}
\left\|\frac{\partial \mathbf{f}(\tau) }{\partial \tau}\bigg\vert_{\tau = t}\right\| \geq \frac{\beta_2^* \beta_3^*}{R_{max}} \left\| \begin{bmatrix} \mathbf{d_2} \\ d_4 \end{bmatrix} \right\| = \alpha_3^* > 0,
\label{eq:uniform_t_d2d4_bound5}
\end{align} for all $t\geq t_0$, with $\beta_2^*$ and $\beta_3^*$ positive constants, and when $\mathbf{d}_3 = \zeros$, $\mathbf{d}_1 = \zeros$ and either $\mathbf{d}_2\neq\zeros$ or $d_4\neq0$. The upper boundedness on the norm of the second derivative of $\mathbf{f}(\tau)$ becomes straightforward under Assumption \ref{assum:bounded_signals}, thus allowing the use of Lemma \ref{lemma:batista_intfunc} in \eqref{eq:uniform_t_d2d4_bound5}. Thus, in this case, it can be shown, using Lemma \ref{lemma:batista_intfunc}, that there exist $\alpha_4^* > 0$ and $\delta_1^* > 0$ such that $\left\|\mathbf{f}(t + \delta_1^*)\right\| \geq \alpha_4^*$ for all $t\geq t_0$. Using Lemma \ref{lemma:batista_intfunc} again, there exist positive constants $\alpha_5^*>0$ and $\delta^*>0$ such that
\begin{equation}
\mathbf{d}^T\mathcal{W}(t,t+\delta^*)\mathbf{d} \geq \alpha_5^*,
\label{eq:uniform_t_d2d4_bound_dWd}
\end{equation} for all $t\geq t_0$ and when $\mathbf{d_1} = \zeros$, $\mathbf{d_3} = \zeros$, and $\mathbf{d_2} \neq \zeros$ or $d_4 \neq 0$. When all the components of $\mathbf{d}$ are null except $d_5$ it follows that 
\begin{equation}
\mathbf{f}(\tau) = \begin{bmatrix}
\mathbf{C_0}\Phi_{BC2}(\tau,t)d_5 \\
\mathbf{C_2}\Phi_{BC2}(\tau,t)d_5 \\
\end{bmatrix},
\label{eq:uniform_fd5_1}
\end{equation} which norm is clearly bounded by
\begin{align}
\|\mathbf{f}(\tau)\| \geq \sigma_{min}(\mathbf{C_0})\left\|\mathbf{\Phi}_{BC2}(\tau,t)d_5 \right\| \label{eq:ftau_bound1}
\end{align} for all $t\geq t_0$. Expanding \eqref{eq:ftau_bound1} and using \eqref{eq:transitionmatrix_1BCm} yields
\begin{align}
\|\mathbf{f}(\tau)\| &\geq \beta_2^* \sqrt{ \sum_{i=1}^{n_r} \left( d_5 \int_{t}^{\tau}{\frac{\sigma-t}{\rhon{1}(\sigma)}d\sigma} \right)^2 }
\label{eq:uniform_fd5_3}
\end{align} for all $t\geq t_0$. In particular at $\tau = t + \delta_2^*$, with $\delta_2^* > 0$, comes
\begin{align}
\|\mathbf{f}(t + \delta_2^*)\| &\geq \beta_2^* \lvert d_5 \rvert \sqrt{ \sum_{i=1}^{n_r} \left( \int_{t}^{t + \delta_2^*}{\frac{\sigma-t}{\rhoi(\sigma)}d\sigma} \right)^2 }
\label{eq:uniform_fd5_4}
\end{align} for all $t\geq t_0$. By the Integral Mean Value theorem there exists $c\in]t,t+\delta_2^*[$ such that 
\begin{equation}
\int_{t}^{t + \delta_2^*}{\frac{\sigma-t}{\rhon{1}(\sigma)}d\sigma} = \delta_2^*\frac{c-t}{\rhoi(c)}.
\label{eq:uniform_fd5_int1}
\end{equation} for all $t \geq t_0$. Now defining $\delta_3^* := c - t > 0$, which is clearly positive because $c > t$, allows to write
\begin{equation}
\int_{t}^{t + \delta_2^*}{\frac{\sigma-t}{\rhon{1}(\sigma)}d\sigma} = \frac{\delta_2^*\delta_3^*}{\rhoi(t+\delta_3^*)}.
\label{eq:uniform_fd5_int2}
\end{equation} for all $t \geq t_0$. Under Assumption \ref{assum:vehicle_motion} comes from \eqref{eq:uniform_fd5_4} and \eqref{eq:uniform_fd5_int2} that
\begin{align}
\|\mathbf{f}(t + \delta_2^*)\| &\geq \beta_2^* \lvert d_5 \rvert \sqrt{ \sum_{i=1}^{n_r} \left( \frac{\delta_2^*\delta_3^*}{R_{max}} \right)^2 } \nonumber \\
&= \frac{\delta_2^*\delta_3^*\beta_2^* \lvert d_5 \rvert \sqrt{ n_r }}{R_{max}} = \alpha_6^* > 0
\label{eq:uniform_fd5_5}
\end{align} for all $t\geq t_0$. Finally, Lemma \ref{lemma:batista_intfunc} is used once again to show that there exist $\alpha > 0$ and $\delta > 0$, for all $t \geq t_0$ and $\{\mathbf{d}\in\mathbb{R}^{8+n_r}:\|\mathbf{d}\| = 1\}$, such that $\mathbf{d}^{\scriptscriptstyle T}\mathcal{W}(t,t+\delta)\mathbf{d} \geq \alpha$, which means that the system is uniformly completely observable and therefore concludes the proof.
\end{proof}
\end{thm}

To recover the augmented system dynamics in the original coordinate space, the original Lyapunov state transformation \eqref{eq:statetransformation1} is reverted considering the augmented state transformation $\mathbf{\Gamma}(t)  := \mathbf{T_r}(t) \mathbf{x}(t)$, where $\mathbf{T_r}(t) := \diag\left(\R^T(t),\R^T(t),1,\ldots,1\right)$ is a also Lyapunov state transformation that preserves all observability properties of the LTV system \eqref{eq:augmenteddynamics_fullmatrixform}.

Including system disturbances and sensor noise yields the final reverted augmented dynamics
\begin{equation*}
\left\{
\begin{aligned}
\mathbf{\dot{\Gamma}}(t) &= \mathbf{A_\Gamma}(t)\mathbf{\Gamma}(t) + \mathbf{B_\Gamma}(t)\velr(t) + \mathbf{n_x}(t), \\
\mathbf{y}(t) &= \mathbf{C_\Gamma}(t)\mathbf{\Gamma}(t) + \mathbf{n_y}(t),
\end{aligned}
\right.
\end{equation*} where
\begin{align*}
\label{eq:sysdyn_matrixA_revert}
&\mathbf{A_\Gamma}(t) = \nonumber \\
&\begin{bmatrix}
-\skews{\w(t)} & -\eye & \zeros & \zeros & \zeros \\
\zeros & -\skews{\w(t)} & \zeros & \zeros & \zeros \\
\frac{\bn{1}^T\skews{\w(t)} - \velr^T(t)}{\rhon{1}(t)} &  \frac{\bn{1}^T}{\rhon{1}(t)} & \zeros & -\frac{1}{\rhon{1}(t)} & 0 \\
\vdots & \vdots & \vdots & \vdots & \vdots \\
\frac{\bn{n_r}^T\skews{\w(t)} - \velr^T(t)}{\rhon{n_r}(t)} &  \frac{\bn{n_r}^T}{\rhon{n_r}(t)} & \zeros & -\frac{1}{\rhon{n_r}(t)} & 0 \\
\zeros & -\velr^T(t) & \zeros & 0 & -1 \\
\zeros & \zeros & \zeros & 0 & 0
\end{bmatrix},
\end{align*}
\begin{equation*}
\label{eq:sysdyn_matrixB_revert}
\mathbf{B_\Gamma}(t) =
\begin{bmatrix}
-\eyen{3} & \zeros & \frac{\bn{1}}{\rhon{1}(t)} & \cdots & \frac{\bn{n_r}}{\rhon{n_r}(t)} & \zeros & \zeros
\end{bmatrix}^T,
\end{equation*}
\begin{equation*}
\mathbf{C_\Gamma}(t) = 
\begin{bmatrix} 
\zerosn{n_r}{3} & \zerosn{n_r}{3} & \mathbf{C_0} & \zerosn{n_r}{2} \\ 
\mathbf{C_1}(t)\R(t) & \zerosn{n_C}{3} & \mathbf{C_2} & \zerosn{n_C}{2}
\end{bmatrix}, 
\end{equation*} where $\mathbf{n_x}(t)$ and $\mathbf{n_y}(t)$ are assumed to be uncorrelated, zero-mean, white Gaussian noise, with $E\left[ \mathbf{n_x}(t) \mathbf{n_x}^T(\tau) \right] = \mathbf{Q_x}(t)\delta(t-\tau)$ and $E\left[ \mathbf{n_y}(t) \mathbf{n_y}^T(\tau) \right] = \mathbf{Q_y}(t)\delta(t-\tau)$.

\section{Simulation Results}
\label{sec:simres}
\begin{figure}[htb]
	\centering
		\includegraphics[width=0.60\columnwidth]{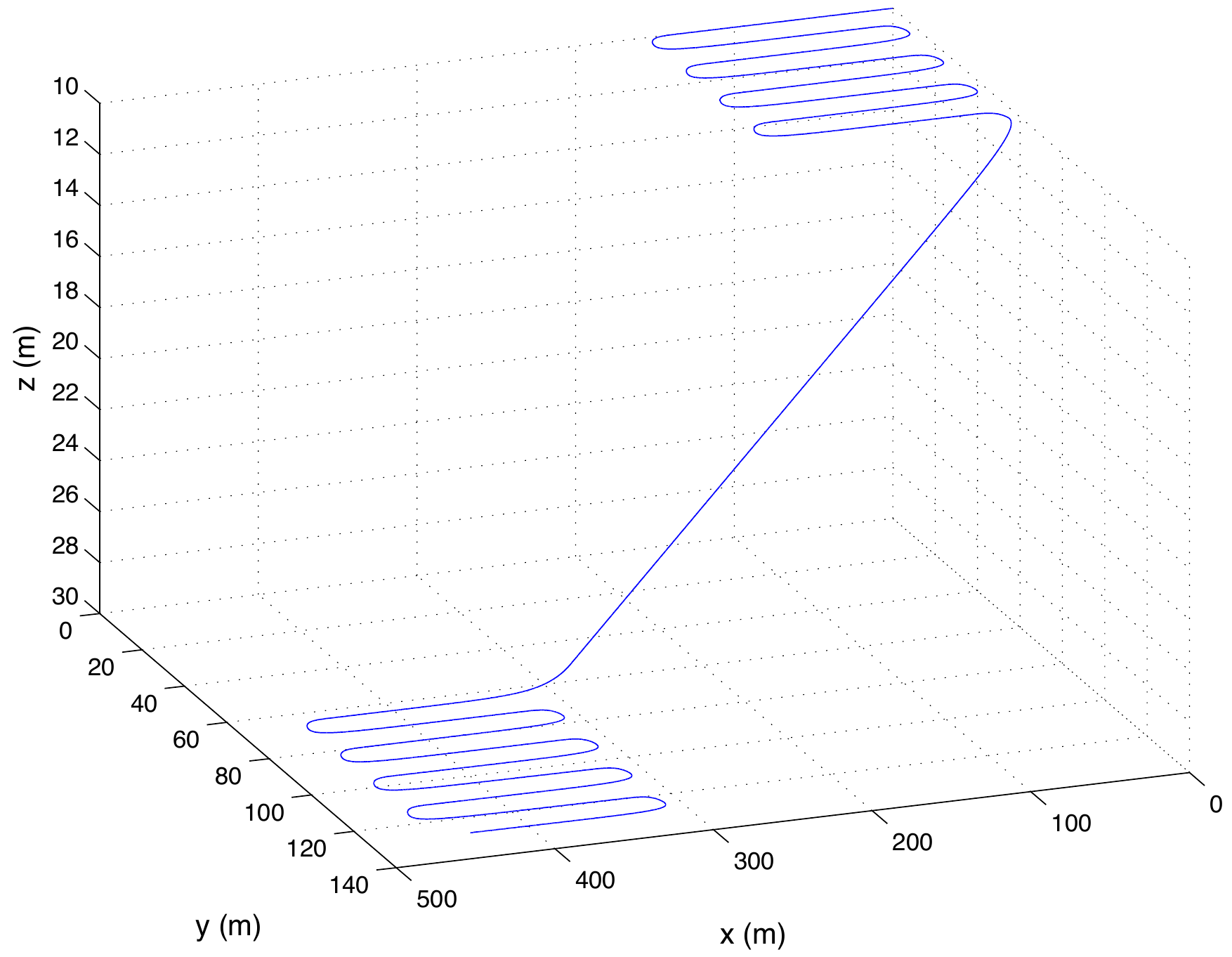}
	\caption{Vehicle nominal trajectory}
	\label{fig:vtraj} 
\end{figure}

The performance of the proposed filter was assessed in simulation using a kinematic model for an underwater vehicle. The vehicle describes a typical survey trajectory as depicted in Figure \ref{fig:vtraj}.

The USBL receiving array is composed of 4 receivers that are installed on the vehicle with an offset of $30$ cm along the $x$-axis of the body-fixed coordinate frame $\BCF$ where the DVL and the rate gyros are also installed. Thus the positions of the receivers with respect $\BCF$ are given in meters by $\bn{1} = \begin{bmatrix} 0.2 & -0.15 & 0  \end{bmatrix}^{\scriptscriptstyle T}$, $\bn{2} = \begin{bmatrix} 0.2 & 0.15 & 0  \end{bmatrix}^{\scriptscriptstyle T}$, $\bn{3} = \begin{bmatrix} 0.4 & 0 & 0.15  \end{bmatrix}^{\scriptscriptstyle T}$, and $\bn{4} = \begin{bmatrix} 0.4 & 0 & -0.15  \end{bmatrix}^{\scriptscriptstyle T}$.

The DVL fluid-relative velocity measurements are considered to be corrupted by additive uncorrelated zero-mean white Gaussian noise with an accuracy of $0.2\%$ of the nominal velocity with an additional standard deviation of $1$ $mm/s$, which is inspired on the LinkQuest NavQuest 600 Micro DVL sensor package. The rate gyros are also inspired on a realistic sensor package, the Silicon Sensing CRS03 triaxial rate gyro, and are thus considered to be disturbed by additive, uncorrelated, zero-mean white Gaussian noise, with a standard deviation of $0.05$ $deg/s$.

The range measurements between the transponder and the reference receiver (receiver $1$) are considered to be disturbed by additive, zero-mean white Gaussian noise, with $1$ $m$ standard deviation whilst the RDOA between receiver $1$ and the other 3 receivers is considered to be measured with an accuracy of $6$ $mm$. The transponder is located in inertial coordinates at ${}^I\mathbf{p}_t = \begin{bmatrix} 200 & 0 & 0 \end{bmatrix}^T$ $[m]$, and the unknown underwater current velocity has an intensity of $0.2$ $m/s$ in all three axis. The augmented states that correspond to the ranges $\xn{3},\ldots,\xn{2+n_r}$ are initialized with the first set of measurements available, the filter position estimate is initialized with an offset of $20$ $m$ from the nominal position, and the rest of the initial estimates is set to zero.

\begin{figure}[htb]
	\centering
		\includegraphics[width=0.90\columnwidth]{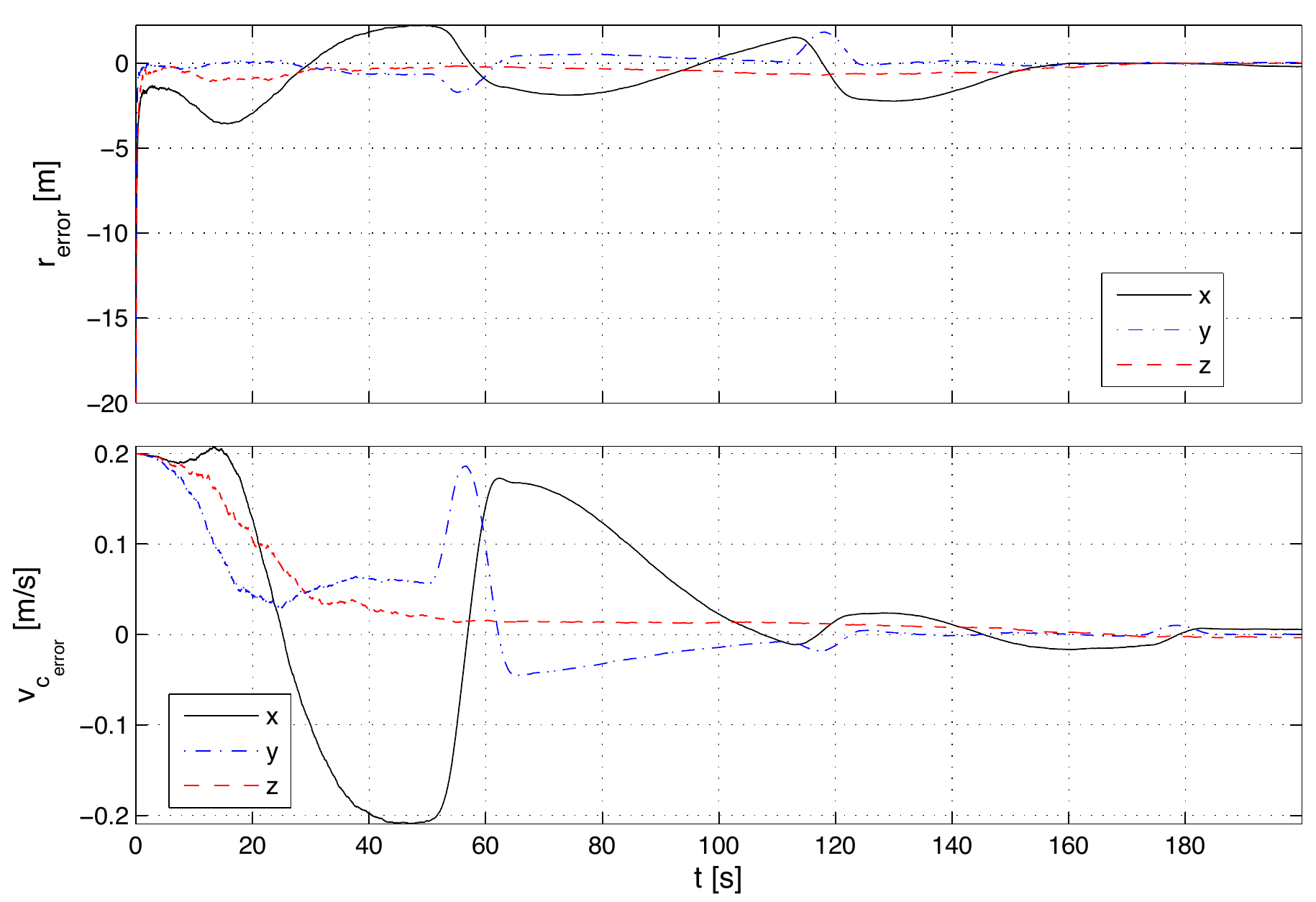}
	\caption{LTV Kalman filter initial convergence}
	\label{fig:ltvkf_est_initial}
\end{figure}

\begin{figure}[htb]
	\centering
		\includegraphics[width=0.90\columnwidth]{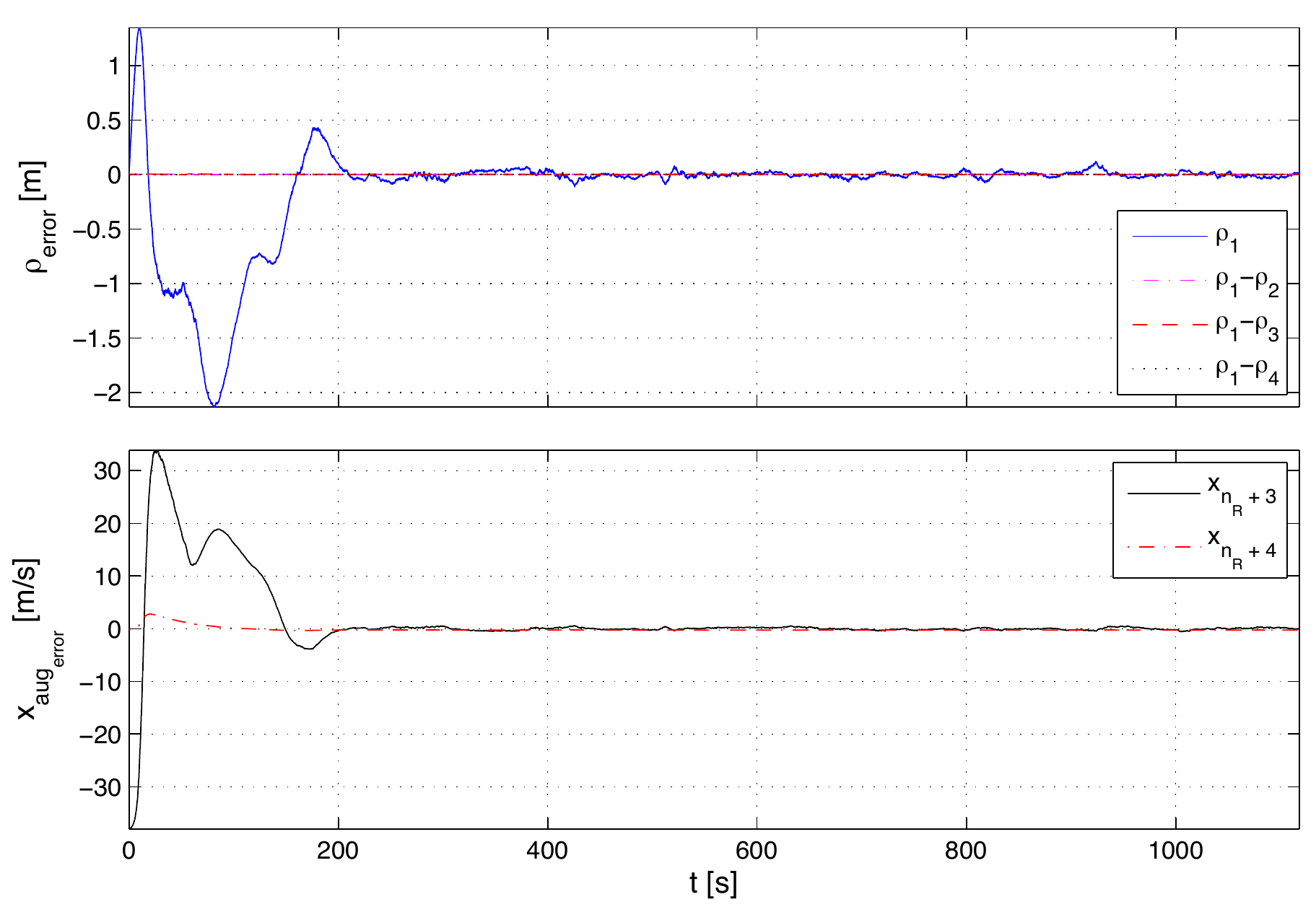}
	\caption{LTV Kalman filter augmented states evolution}
	\label{fig:ltvkf_estaug_total} 
\end{figure}

\begin{figure}[htb]
	\centering
		\includegraphics[width=1.0\columnwidth]{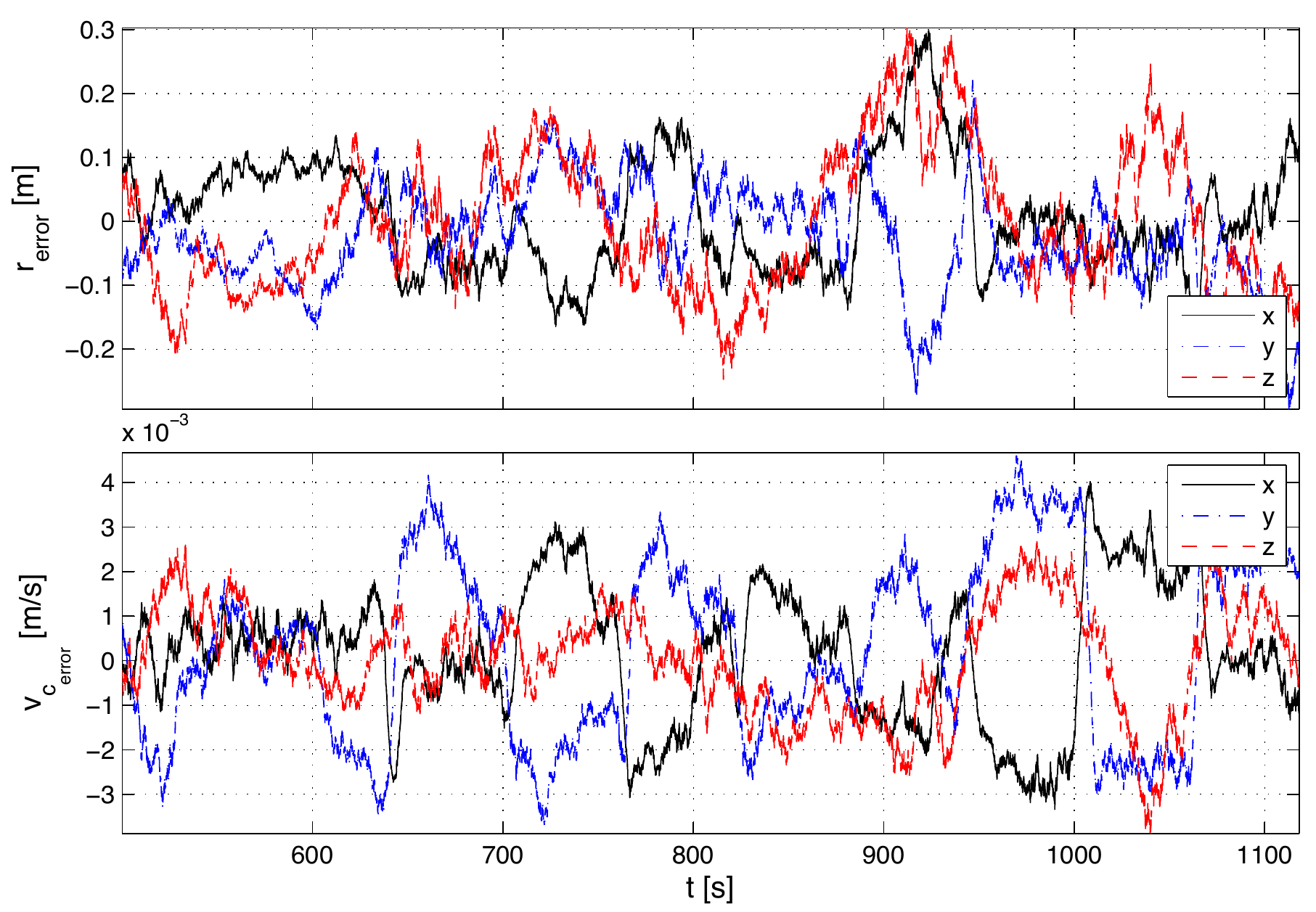}
	\caption{LTV Kalman filter steady state response}
	\label{fig:ltvkf_est_steadystate} 
\end{figure}

The initial evolution of the position and current velocity estimation error is depicted in Figure \ref{fig:ltvkf_est_initial}. The full evolution of the augmented states error is represented in Figure \ref{fig:ltvkf_estaug_total}. The attainable performance of the filter is better illustrated in Figure \ref{fig:ltvkf_est_steadystate}, where the steady-state response of the position and current velocity estimation error is shown.

The performance of the filter is compared with two alternative filter designs, the well-known and established Extended Kalman Filter (EKF) and the Kalman filter with the Planar-Wave approximation (KFPW). The first design linearizes the nonlinear range and RDOA measurements about the filter estimates in order to compute a suboptimal Kalman gain. In the latter, the feedback is accomplished by means of a precomputed transponder position fix from the USBL that resorts to a planar-wave approximation, previously used by the authors \cite{Fusion06}.

The Root-Mean-Square (RMS) of the steady-state position error of the three filters is presented in Table \ref{tab:steadystate_rms_comp}. Comparing the steady-state response of the three filters it can be seen that all of them attain the same performance level. The solution presented in this work has the advantage of being GAS, which is not guaranteed for the other designs.

\begin{table}[htb]
	\centering
	\caption{Steady-state transponder position $\rpb(t)$ error RMS comparison}
	\begin{tabular}{|r|c|c|c|}
			\hline
		  Filter & $x$ $[m]$ & $y$ $[m]$ & $z$ $[m]$ \\ \hline
			LTV Kalman & 0.0880 & 0.0837 &  0.1101 \\ \hline
			KFPW & 0.1122 & 0.0847 & 0.1096 \\ \hline
			EKF & 0.0611  & 0.0814 & 0.1128 \\ \hline
		\end{tabular}
	\label{tab:steadystate_rms_comp}
\end{table}

\section{Conclusions}
\label{sec:conclusions}
The main contribution of the paper lies on the design of globally asymptotically stable position filters based directly on the nonlinear sensor readings of USBL and a DVL. At the core of the proposed filtering solution is the derivation of a LTV system that fully captures the dynamics of the nonlinear system. The LTV model is achieved through appropriate state augmentation allowing for the use of powerful linear system analysis and filtering design tools that yield GAS filter error dynamics.

The performance of the proposed filter was assessed in simulation and compared against two traditional solutions, the EKF and a Kalman filter that resorts to the planar approximation of the acoustic wave arriving at the USBL array. Comparison of the steady-state position error from the three designs lead to the conclusion that all demonstrated similar performance levels using realistic sensor noise and disturbances. The advantage of the new filter structure is nevertheless evident, due to its GAS properties which is not guaranteed for either of the two traditional solutions.

\bibliographystyle{IEEEtran}
\bibliography{cdc2010_extendedversion}

\end{document}